\documentclass{amsart}
\newtheorem{theorem}{Theorem}[section]
\newtheorem{lemma}[theorem]{Lemma}
\theoremstyle{definition}
\newtheorem{definition}[theorem]{Definition}

\theoremstyle{remark}
\newtheorem{remark}[theorem]{Remark}

\numberwithin{equation}{section}

\begin{document}
\def\C{\mathbb C}
\def\D{\mathbb D}
\def\R{\mathbb R}
\def\X{\mathbb X}
\def\cA{\mathcal A}
\def\cT{\mathcal T}
\def\Z{\mathbb Z}
\def\Y{\mathbb Y}
\def\Z{\mathbb Z}
\def\N{\mathbb N_0}
\def\cal{\mathcal}
\def\F{\mathcal F}


\title[Traveling waves]{Traveling wave dispersal in partially sedentary age-structured populations}
\author{Thuc Manh  Le }
\address{Department of Mathematics, Vietnam National University, 334 Nguyen Trai, Hanoi, Vietnam}
\email{manhthucle@gmail.com}
\author{Frithjof Lutscher}
\address{Department of Mathematics and Statistics, University of Ottawa,
585 King Edward Avenue, Ottawa, ON K1N 6N5, Canada}
\email{ flutsche@uottawa.ca}
\author{Nguyen Van Minh}
\address{Department of Mathematics, University of West Georgia, Carrollton, GA 30118, USA}
\email{vnguyen@westga.edu}

\date{\today}
\begin{abstract}
In this paper we present a thorough study on the existence of traveling waves in a mathematical model of dispersal in a partially sedentary age-structured population. This type of model was first proposed by Veit and Lewis in [{\it Am. Nat.}, {\bf 148} (1996), 255-274]. We choose the fecundity function to be the Beverton-Holt type function. We extend the theory of traveling waves in the population genetics model  of Weinberger in [{\it SIAM J. Math. Anal.}, {\bf 13} (1982), 353-396] to the case when migration depends on age groups and a fraction of the population does not migrate.
\end{abstract}
\keywords{Traveling waves; spreading speed; partially sedentary population; delay effect}
\subjclass[2010]{Primary: 92D25; Secondary: 37N25; 39A22}

\maketitle


\section{Introduction}
In \cite{veilew} a mathematical model for dispersal in a partially, sedentary age-structured population
was developed to simulate the spatial spread of the house finch ({\em carpodacus mexicanus}). The house finch is native to the southern part of the United States and to Mexico. It spread quickly in the eastern part of the United States and Canada in the 1940s after the release of captive specimens in the New York City area. The model is of the form
\begin{eqnarray}\label{eq1}
N_{n+1} (x) &=& s(1-p_A)N_n (x) +(1-p_J) F(N_n(x))+\int^\infty_{-\infty} K_A(|x-y|) sp_AN_n(y)dy \nonumber\\
&& + \int^\infty_{-\infty} K_J(|x-y|) p_J F(N_n (y))dy , \quad n=1,2,\cdots ,
\end{eqnarray}
where $N_n =J_n+A_n$ is the sum of the juvenile and adult bird densities in year $n$. It is assumed that juvenile birds born in one season mature by the next season. Adults survive from one season to the next with probability $s.$ The number of surviving offspring born in year $n$ is denoted by the fecundity function $F(N_n).$ Juveniles and adults may differ in their probability to disperse as well as in their dispersal behavior if they disperse. The fractions of dispersing juveniles and adults are denoted by  $p_J$ and $p_A,$ respectively. The dispersal behavior is given by the probability density functions $K_J, K_A,$ also called  dispersal kernels.

In this paper, we will always assume that the fecundity function is a monotone increasing and bounded function, for example a function of Beverton-Holt type or, equivalently, of Holling type II. The Beverton-Holt function originated as a stock-recruitment function in fisheries and is now very common in classic discrete-time population models. By appropriate rescaling, we can always assume that $F$ is of the form
 \begin{equation}\label{1.2}
F(u):=\frac{krMu}{M+(r-1)u},
\end{equation}
with $k+s=1.$ All parameters are assumed positive. At low population density, the per capita number of offspring is $kr,$ and the carrying capacity of the non-spatial model $N_{n+1}=sN_n+F(N_n)$ is $M.$

\bigskip
As noted in \cite{vollui}, although this model was used to simulate the migration of house finches in \cite{veilew}, no mathematical analysis was given yet. A first attempt to study  models with partially sedentary populations was made in \cite{vollui}. These authors studied the special case that dispersal probability and dispersal behavior are independent of age-structure. More precisely, under the assumptions $p_J=p_A$ and $K_J=K_A,$ the above model falls within the framework considered in  \cite{vollui}. It is the purpose of this paper to give a thorough study of the asymptotic behavior of solutions of Eq.~(\ref{eq1}), (\ref{1.2}) with $p_A\ne p_J$ and $K_J\ne K_A.$ Qur goal is to extend the fundamental theory for spreading speeds and traveling waves developed by Weinberger to this case. For the existence of traveling waves the fundamental assumption on compactness of \cite[Theorem 6.6]{wei} is not satisfied for Eq. (\ref{eq1}), nor the weak compactness condition in \cite{liazha2}. We will use a different approach to prove the existence of traveling waves in this case. This idea is also employed in \cite{lutmin}. Our main result is Theorem \ref{the wave} that complements some results in \cite{vollui}.

To simplify notation, we will assume that $K_J=K_A=\colon K$ from here on. We return to the general case in the discussion.

\subsection*{Notations and Assumptions}
We denote by $\R$ the real line. We also denote by $BM(\R ,\R)$ ($BC(\R,\R)$, respectively) the space of all measurable and essentially bounded real valued functions on $\R$ (the space of all bounded continuous real valued functions on $\R$, respectively) with essential sup-norm. For a constant $\alpha$ we will denote the constant function $\R \ni x \mapsto \alpha$ by this number $\alpha$ for convenience if this does not cause any confusion.  $C_M$ stands for the set $\{ f\in BC(\R,\R)| f(x)\in [0,M]\}$, and $BM(\R,[0,M]):= \{ f\in BM(\R,\R)| f(x)\in [0,M]\}$. The metric on  $C_M$ is defined by the sup norm. In $BM(\R,\R)$ we use the natural order defined as  $u\le v$ if and only if $u(x)\le v(x)$ for all $x\in \R$.

\bigskip
Unless otherwise stated, we assume that the parameters in the function $F$  satisfy $r>1, M>0.$
We also assume that $K(|x|)$ is a probability density function defined on $\R$ and satisfies
\begin{equation}
\int^\infty_{-\infty} e^{\mu x}K(|x|)dx <\infty , \quad \mbox{for all} \ \mu\in \R .
\end{equation}

\section{Spreading Speed}
Most of the results in this section are derived from the general theory on spreading speeds in \cite{wei}. For later use in the paper we will discuss details of these results below. Let us define a dynamical system $u_{n+1}=Q[u_n]$  by setting
\begin{eqnarray}\label{eq2}
Q[u](x) &=&s(1-p_A)u (x) +(1-p_J) F(u(x))+\int^\infty_{-\infty} K(|x-y|) sp_Au(y)dy \nonumber\\
&& + \int^\infty_{-\infty} K(|x-y|) p_J F(u (y))dy ,
\end{eqnarray}
for each $u\in BM(\R,\R)$.
\begin{lemma}\label{lemma 2.1}
Under the above notations and assumptions, the operator $Q$ is an operator acting in $BM(\R,\R)$ leaving $BC(\R,\R)$ invariant with the following properties
\begin{enumerate}
\item $Q[0]=0$, $Q[M]=M$, $Q[\alpha ]>\alpha $ for all $0<\alpha <M$;
\item If $u,v\in BM(\R,\R)$ such that $u\ge v$, then $Q[u]\ge Q[v]$;
\item If $u_n\in BC(\R,\R)$ such that $u_n$ is convergent to $u$ uniformly on each bounded subset of $\R$, then $Q[u_n](x)$ is convergent to $Q[u](x)$ for each $x\in \R$;
\item If $\alpha >M$, then $Q[\alpha ] < \alpha$;
\item There is a constant $\bar \gamma$ such that $\gamma < Q[\gamma ] < \bar \gamma$ for all $\gamma \in [0,\bar \gamma )$, and $Q[\bar\gamma ] =\bar \gamma$.
\end{enumerate}
\end{lemma}
\begin{proof} Before proving the properties of $Q$ we notice that from the definition of $Q$ it maps $BM(\R,\R)$ into itself. Next, to show that it leaves $BC(\R,\R)$ invariant it is sufficient to prove that the integrals
\begin{eqnarray}\label{eq103}
\int^\infty_{-\infty} K(|x-y|)  u(y) dy &=& -\int^\infty_{-\infty} K(|\xi |)   u(x-\xi ) d\xi \\
\int^\infty_{-\infty} K(|x-y|)  F(u(y)) dy &=& -\int^\infty_{-\infty} K(|\xi |)   F(u(x-\xi )) d\xi
\end{eqnarray}
depend continuously on $x$. To this purpose, for each $\epsilon>0$ we can find a large enough $T>0$ such that
\begin{eqnarray}\label{eq104}
\int^{-T}_{-\infty} K(|\xi |)d\xi   + \int^\infty_{T} K(|\xi |)d\xi  <\frac{\epsilon }{4 (\| u\|+\| F(u)\| )}.
\end{eqnarray}
On the interval $[-T-1,T+1]$ the functions $u$ and $F(u)$ are both uniformly continuous. Therefore, for a given $x_0$ there exists $1>\delta>0$ such that if $|x_0-x|<\delta$, then, for all $\xi\in [-T,T]$
\begin{eqnarray}
 | u(x_0-\xi )-   u(x_0-\xi )|<\frac{\epsilon }{4}\\
  |F( u(x_0-\xi ))-   F(u(x_0-\xi ))|<\frac{\epsilon }{4}.
\end{eqnarray}
This yields that for all $|x_0-x|<\delta$,
\begin{eqnarray}
&&| \int^\infty_{-\infty} K(|\xi |)   u(x_0-\xi ) d\xi -\int^\infty_{-\infty} K(|\xi |)   u(x-\xi ) d\xi | \nonumber \\
&& \le |\int^\infty_{T} K(|\xi |)    d\xi +\int^\infty_{T} K(|\xi |)   d\xi |\cdot  \| u\| \nonumber \\
&&+ \int^T_{-T} K(|\xi |)  | u(x_0-\xi ) - u(x-\xi )| d\xi |  \nonumber  \\
&& \le \frac{\epsilon}{4}+ \frac{\epsilon}{4}  < \epsilon . \nonumber
\end{eqnarray}
Similarly
\begin{eqnarray}
&&| \int^\infty_{-\infty} K(|\xi |)   F(u(x_0-\xi )) d\xi -\int^\infty_{-\infty} K(|\xi |)   F(u(x-\xi )) d\xi |  < \epsilon\nonumber .
\end{eqnarray}
That is the continuity of the integrals in (\ref{eq103}), proving that $Q$ leaves $BC(\R,\R)$ invariant.

\bigskip
Now we prove the properties of $Q$:\\
(i): This property is clear because of the assumption $s+k=1.$ Indeed, $Q[0]=0$, and
\begin{eqnarray}\label{eq100}
Q[M](x) &=&s(1-p_A)M +(1-p_J) \frac{krMM}{M+(r-1)M}+\int^\infty_{-\infty} K(|x-y|) sp_AMdy \nonumber\\
&& + \int^\infty_{-\infty} K(|x-y|) p_J \frac{krMM}{M+(r-1)M}dy \nonumber\\
&=& s(1-p_A)M + (1-p_J) kM +sp_AM+p_JkM   \nonumber     \\
&=& (s+k)M \nonumber    \\
&=& M.
\end{eqnarray}
Notice that $F(\alpha ) >k\alpha$ for all $0<\alpha <M$. Therefore, for all $0<\alpha <M$
\begin{eqnarray}\label{eq100}
Q[\alpha ](x) &>&s(1-p_A)\alpha +(1-p_J) k \alpha +\int^\infty_{-\infty} K(|x-y|) sp_A\alpha
dy \nonumber\\
&& + \int^\infty_{-\infty} K(|x-y|) p_J k \alpha
dy \nonumber\\
&=& s(1-p_A)\alpha + (1-p_J) k\alpha +sp_A\alpha +p_Jk\alpha    \nonumber     \\
&=& (s+k)\alpha \nonumber    \\
&=& \alpha .
\end{eqnarray}
(ii): This property follows since the function $F$ is increasing.\\
(iii): This property is clear from the definition of $Q$.\\
(iv): We have
\begin{eqnarray}\label{eq101}
F'(u)-k &=& \frac{kr(M+(r-1)u)-kru(r-1)}{(M+(r-1)u)^2} -k\nonumber\\
&=& \frac{krM}{(M+(r-1)u)^2} -k  \nonumber.
\end{eqnarray}
Notice that $F'(u)-k <0$ if $u>M$. Therefore, $F(\alpha )< k \alpha $ if $\alpha >M$. This yields that if $\alpha >M$
\begin{eqnarray}\label{eq102}
Q[\alpha ](x) &<&s(1-p_A)\alpha +(1-p_J) k \alpha +\int^\infty_{-\infty} K(|x-y|) sp_A\alpha
dy \nonumber\\
&& + \int^\infty_{-\infty} K(|x-y|) p_J k \alpha
dy \nonumber\\
&=& \alpha .
\end{eqnarray}
(v): Actually we can choose $\bar\gamma =M$. In fact, for each $\gamma \in (0,M)$ we have $F(\gamma )>k\gamma$. Therefore, by the above computation (\ref{eq100}), $\gamma<Q(\gamma )<Q[M]=M$ for all $\gamma \in [0,M]$.
\end{proof}

Next,  we apply the general theory on spreading speeds in \cite{wei} to our dynamical system
\begin{equation}
u_{n+1}=Q[u_n], \quad n=1,2,\cdots .
\end{equation}
Basically, according to the theory of \cite{wei}, for our model the following procedure of defining the spreading speed is valid: let us choose a function $\varphi :\R \to \R$ such that
\begin{enumerate}
\item $\varphi$ is continuous and non-increasing,
\item $\varphi (-\infty ):= \lim_{t\to-\infty}\varphi (t) \in (0,M)$,
\item $\varphi (s)=0$ for $s\ge 0$.
\end{enumerate}
We then define an operator $R_c[\cdot ]$ on the space $C_M$ for every constant $c$ as
\begin{equation}
R_c[u] (s) := \max \{ \varphi (s), Q[u(c+\cdot) ](s)\}, \quad s\in \R ,
\end{equation}
and a sequence of functions $\{ a_n (c;\cdot )\}$ by
\begin{equation}
a_{n+1}:= R_c[a_n],\quad a_0=\varphi .
\end{equation}
As shown in \cite{wei} the sequence $\{ a_n (c;\cdot )\}$ is increasing and bounded, so for each $s\in\R$, we obtain the pointwise limit
\begin{equation}\label{def a}
\lim_{n\to\infty} a_n(c;s)=a(c;s), \quad s\in \R.
\end{equation}
Obviously, $0\le a(c;s) \le M$ for all $s\in \R$. The following number is called {\it spreading speed} for our model
\begin{equation}
c^*:= \sup \{ c| a(c;+\infty )=M\} .
\end{equation}
Applying the theory in \cite[Lemma 5.2 and Proposition 5.1]{wei} to our model gives the following:
\begin{lemma} \label{lem boundary cond}
\begin{enumerate}
\item For each $c\in \R$
\begin{equation}
a(c;-\infty )=M;
\end{equation}
\item If $c\ge c^*$, then
\begin{equation}
a(c;+\infty )=0.
\end{equation}
\end{enumerate}
\end{lemma}
Weinberger \cite{wei} proved that if there exists a bounded non-negative measure $m(x,dx)$ on $\R$ such that
\begin{equation}\label{2.8}
Q[u](x) \le \int^\infty_{-\infty} u(x-y)m(y, \ dy ), \quad u\in C_M,
\end{equation}
then
\begin{equation}
c^* \le \inf_{\mu >0} \frac{1}{\mu } \ln \int^\infty_{-\infty} e^{\mu x}m(x,\ dx) .
\end{equation}
And if there exists a bounded non-negative measure $l(x,\ dx)$ with property that $\int^\infty_{-\infty} l(x,\ dx) >1$ and $Q[u](x) \ge \int u(x-y)l(y, \ dy)$ for all $u$ such that $0\le u(x)\le \epsilon$, then
\begin{equation}\label{2.9}
c^* \ge \inf_{\mu >0} \frac{1}{\mu} \ln \int^\infty_{-\infty} e^{\mu x}l(x,\ dx).
\end{equation}

Below we follow the argument of \cite{vollui} to give an estimate of the spreading speed. We see that since for $u\in C_M$ $F(u) \le kru$,
\begin{equation}
Q[u](x)  \le \left[ s(1-p_A)  +(1-p_J) kr \right]u(x)+\int^\infty_{-\infty} K(|x-y|) \left( sp_A +p_J kr  \right) u(y)dy.
\end{equation}
If we let
$$
m(x,\ dx )= \left[ s(1-p_A)  +(1-p_J) kr\right]\delta_0+K(|x|) \left( sp_A +p_J kr \right) ,
$$
where $\delta_0$ is the Dirac delta measure, then (\ref{2.8}) holds. On the other hand, for each $1<r_1<r$ there exists $\epsilon >0$ such that $F(u)\ge kr_1 u$ for $u$ such that $0\le u(x)\le \epsilon$ for all $x\in\R$. Therefore, for such $u$
\begin{equation}
Q[u](x)  \ge \left[ s(1-p_A)  +(1-p_J) r_1\right]u(x)+\int^\infty_{-\infty} K(|x-y|) \left( sp_A +p_J r_1 \right) u(y)dy.
\end{equation}
We let
$$
l(x, \ dx):= \left[ s(1-p_A)  +(1-p_J) kr_1\right]\delta_0+ K(|x|)  ( sp_A +p_J kr_1 ).
$$
Then, $\int_{-\infty}^\infty l(x,\ dx)>1$ and (\ref{2.9}) holds. Next, since $r_1$ can be chosen arbitrarily between $1$ and $r$, we have
\begin{equation}\label{sp sp}
c^*= \inf_{\mu>0}\frac{1}{\mu} \int^\infty_{-\infty}  ( sp_A +p_J kr  )e^{\mu x}K(|x|)dx  + [ s(1-p_A)  +(1-p_J) kr ] .
\end{equation}

\section{Traveling Waves}
This section contains our main results of the paper.
\begin{definition}
A monotone traveling wave solution with speed $c$ connecting $0$ to $M$, or for short a traveling wave, of Eq. (\ref{eq1}) is defined to be a non-increasing continuous function $w$ such that $\lim_{x\to \infty}w(x)=0$, $\lim_{x\to-\infty }w(x)=M$, and
$N_n(x):= w(x-nc)$ is a solution of Eq. (\ref{eq1}).
\end{definition}
If we substitute $N_n$ into (\ref{eq1}), we will have
\begin{eqnarray}
w(x-(n+1) c) &=& s(1-p_A)w(x-nc) +(1-p_J) F(w(x-nc))\nonumber \\
&&+\int^\infty_{-\infty} K(|x-y|) sp_A w(y-nc)dy \nonumber\\
&& + \int^\infty_{-\infty} K(|x-y|) p_J F(w(y-nc))dy , \quad n=1,2,\cdots .
\end{eqnarray}
If we set $\xi := x-(n+1)c$ and $z:=y-nc$ then, the above equation becomes what is called "wave equation" associated with (\ref{eq1}):
\begin{eqnarray}
w(\xi ) &=& s(1-p_A)w(\xi +c) +(1-p_J) F(w(\xi +c))+\int^\infty_{-\infty} K(|\xi +c-z|) sp_A w(z)dz \nonumber\\
&& + \int^\infty_{-\infty} K(|\xi +c-z|) p_J F(w(z))dz , \quad n=1,2,\cdots .
\end{eqnarray}

\bigskip
For each $u\in C_M$, or more generally $BM(\R,\R)$, and $c\in \R$ we set
\begin{eqnarray}\label{eq31}
B_c[u] (x) &=&s(1-p_A)u (x+c) +(1-p_J) F(u(x+c)) \\
C_c[u](x) &=& \int^\infty_{-\infty} K(|x +c-y|) sp_A u(y)dy  \nonumber\\
&&+ \int^\infty_{-\infty} K(|x +c-y|) p_J F(u(y))dy .
\end{eqnarray}
Obviously, a monotone traveling wave solution to (\ref{eq1}) is a non-increasing continuous function $w$ with $w(-\infty)=M, w(\infty )=0$ that is a fixed point of $Q_c:= B_c +C_c$. We remark that although there are several extensions of Weinberger's theory on the existence of monotone traveling waves, to our best knowledge the existence of monotone traveling wave to (\ref{eq1}) is still open. The reason is the operator $Q_c$ does not satisfy any compactness conditions listed in \cite{wei}, \cite{liazha2}. Volkov and Lui extended Weinberger's theory to a class of systems without compactness conditions. However, the model considered in \cite{vollui} includes our model (\ref{eq1}) only if $p_A=p_J$, that means that age structure  does not affect migration behavior. However, it is more realistic to make the assumption that $p_J\not= p_A$. This assumption makes the problem of studying the existence of monotone traveling waves much harder. Below we will prove the existence of monotone traveling waves of (\ref{eq1}) under assumption that the function $F$ in (\ref{eq1}) is of the form (\ref{1.2}).

\begin{lemma}\label{lem b}
 Assume that
\begin{eqnarray}
s(1-p_A)+(1-p_J)kr &<& 1 .\label{b2}
\end{eqnarray}
Then, for each given $w\in C_M$ ($w\in BM(\R,[0,M])$, respectively) the equation operator $u-B_c[u ]=w$ has an unique solution $u$ in $C_M$ (in $BM(\R,[0,M])$, respectively) which will be denoted by $u:=G_cw$.
\end{lemma}
\begin{proof}
Consider the function
\begin{equation}
g(x)= s(1-p_A)x+(1-p_J)\frac{krMx}{M+(r-1)x}, \quad x\in [0,M].
\end{equation}
For all $x\in[0,M]$,
\begin{eqnarray*}
g'(x)&=& s(1-p_A)+ (1-p_J)\frac{krM(M+(r-1)x)-krMx(r-1) }{(M+(r-1)x)^2}\\
&=& s(1-p_A)+ (1-p_J) \frac{krM^2}{(M+(r-1)x)^2}\\
&\le & s(1-p_A)+ (1-p_J) \frac{krM^2}{(M+(r-1)\cdot 0)^2}\\
&=& s(1-p_A)+  (1-p_J) kr\\
&<& 1.
\end{eqnarray*}
Therefore, there exists a positive $0<p<1$ such that
\begin{eqnarray}
 0< \sup_{x\in [0,M]} |g'(x)| <p .
\end{eqnarray}
Next, we solve the equation $u-B_c[u]=w$ for each given $w\in C[0,M]$. Note that in this case, $B_c$ is a strict contraction because
\begin{eqnarray}
\| B_c[u_1]-B_c[u_2]\| &=& \sup_{x\in \R} | g(u_1(x+c))-g(u_2(x+c))| \nonumber\\
&\le& \sup_{\xi \in \R} |g'(\xi )| \cdot |u_1(x+c)-u_2(x+c)| \nonumber\\
&=&p \| u_1-u_2\| .
\end{eqnarray}
Therefore, by a standard argument we can prove the existence of $(I-B_c)^{-1}$ that is Lipschitz continuous.
\end{proof}

\begin{lemma}\label{lem continuity}
Assume that
the kernel $K(|x|)$ satisfies the above mentioned conditions. Then, for each monotonous $u$ the function $C_c[u]$ is continuous.
\end{lemma}
\begin{proof}
The proof can be done in the same way as in that of \cite[Lemma 3.5]{lutmin}.
\end{proof}

\begin{remark}
We notice that although the operator $B_c$ is a strict contraction in the uniform convergent topology it is not a strict contraction in the norm of compact open topology as defined in \cite{liazha2}. Moreover, operator $C-c$ is not compact the norm of compact open topology because the kernel $K$ may not be continuous, so condition (A3) in \cite{liazha2} is not satisfied with the operator $B_c+C_c$. That is, the theory of traveling waves in \cite{wei} as well as its extension in \cite{liazha2} does not apply to this case.
\end{remark}

The following is the main result of the paper:
\begin{theorem}\label{the wave}
Let all assumptions in Lemmas \ref{lem b} and \ref{lem continuity} be satisfied. Then, if $c\ge c^*$, then there exists a monotone traveling wave to Eq. (\ref{eq1}).
\end{theorem}
\begin{proof}
Let the function $a(c;\cdot )$ be defined as in (\ref{def a}).
Set $\phi_1(s):= a(c;s)$. Note that since $a(c;\cdot)$ is non-increasing and bounded it is a measurable and  bounded function on $\R$.
We define a sequence
\begin{eqnarray}
\phi_{n+1}&=&Q_c[\phi_n],\ n\ge 1 , \quad n=1,2,\cdots .
\end{eqnarray}
We now show that $\{ \phi_n\}$ is a non-increasing sequence in $BM(\R,[0,M])$. In fact, by definition of the sequence $\{ a_n(c;\cdot )\}$ we have
\begin{eqnarray}
  a_{n+1}(c;s)&:=& \max \{\varphi^i(s), Q[a_n(c; \cdot +s+c)](0)\}\nonumber \\
  &=& \max \{ \varphi^i(s), Q[a_n(c; \cdot +c)](s)\}  \\
  &\ge &  Q[a_n(c; \cdot +c)](s) ,\nonumber\\
  &=& Q_c[a_n(c; \cdot )](s)
\end{eqnarray}
Therefore,
\begin{eqnarray} \label{6.7}
  a  (c;s)
  &\ge &  Q_c[a (c; \cdot )](s) .
\end{eqnarray}
That is,
\begin{equation}\label{6.8}
\phi_1 \ge \phi_2 .
\end{equation}
 Since $Q_c$ is order-preserving, by the definition of $\phi_{n+1}$, (\ref{6.8}) yields that
 \begin{equation}\label{6.9}
\phi_{n} \ge \phi_{n+1}, \quad n\in\N .
\end{equation}
Therefore, the sequence $\{ \phi_n\}$ is pointwise non-increasing and bounded below by zero (because these functions are non-negative), so it has a limit $W$ that is a non-increasing function, so it is measurable.
We will show that $W$ is a traveling wave solution to Eq. (\ref{eq1}). By Lemma \ref{lem boundary cond}, $W(-\infty )= { M}$. Next, ${ 0} \le W(+\infty )\le \phi_1(+\infty )={ 0}$, so, $W(+\infty )={ 0}$.  In particular, $W$ is a fixed point of $Q_{c}$, that is
\begin{equation}\label{100}
W=Q_c[W]:=B_c[W]+ C_c[W].
\end{equation}
We need only to show that $W$ is continuous for it to be a traveling wave as in our definition. Since (\ref{100}) is equivalent to the following
\begin{equation}\label{101}
W -B_c[W]= C_c[W],
\end{equation}
by Lemma \ref{lem b} it is equivalent to
\begin{equation}\label{102}
W = [I-B_c]^{-1}C_c[W],
\end{equation}
By Lemma \ref{lem continuity} the function $C_c[W]$ is continuous. In turn, by Lemma \ref{lem b}, $G_c[W]=[I-B_c]^{-1}C_c[W]$ is continuous. Therefore, $W$ is continuous, and thus it is a wave solution of Eq.~(\ref{eq1}).
\end{proof}
\begin{remark}
As in \cite{wei,liazha2,lutmin} by using Lemma \ref{lem boundary cond} we can easily show  that if $c<c^*$ the traveling waves do not exist. With this said the spreading speed $c^*$ is exactly the minimal wave speed of traveling wave solution to Eq. (\ref{eq1}).
\end{remark}

\section{Discussion}
In Section 3 we considered the existence of traveling waves when $F$ is of the form (\ref{1.2}). The function $F$ can chosen to be a more general one that satisfies the following conditions with given positive numbers $r,k, M$:
\begin{enumerate}
\item[(H1)] $F\in C^1[0,M]$;
\item[(H2)] $F(0)=0$, and $F(M)=kM$, $s+k=1$;
\item[(H3)] $F(u)>ku$, for $u\in (0,M)$;
\item[(H4)] $F'(u)\ge 0$, and $F'(0)=kr$, where $r>1$ is a given constant;
\item[(H5)] $F(u)\le kru$, for $u\in [0,M]$;
\item[(H6)] $F'(x)$ is non-increasing on $[0,M]$.
\end{enumerate}
Then, the statement of the main results as well as its proof are unchaged.

\bigskip
Results of the previous sections can be easily extended to the case the habitat is multiple dimensional $\R^d$ with $d=2,3,\dots$. There are no big changes in the statements of the results. And the ideas of proofs remain similar.

\bigskip
We turn to the more biological aspects of our work. Condition (\ref{b2}) imposes an upper bound on the per capita offspring production or, more precisely, on the proportion of individuals who do not disperse.
The original model by Veit and Lewis \cite{veilew} allowed for different dispersal behavior of juveniles and adults. The theory presented here easily extends to the case $K_A\ne K_J.$ Most importantly,  Lemma \ref{lem continuity} holds if the continuity conditions holds for both kernels. Formula (\ref{sp sp}) for the spreading speed becomes $c^*= \inf_{\mu>0}\frac{1}{\mu}\kappa(\mu),$ where
\begin{equation}
\kappa(\mu) = \int^\infty_{-\infty}  [sp_A K_A(|x|)+p_J kr K_J(|x|)] e^{\mu x} dx  + [ s(1-p_A)  +(1-p_J) kr ] .
\end{equation}
More important differences between the model by Veit and Lewis and our analysis here is that they considered the function $F$ to describe a strong Allee effect and  the dispersal probabilities, $p_j, p_A,$ to depend on population density. A strong Allee effect occurs if the per capita population growth rate is highest for intermediate population densities so that the population actually declines for small densities. While the existence of a spreading speed in the presence of an Allee effect is still guaranteed by Weinberger's theory, the existence of traveling waves for that case is a wide open question. Equally open is the question of traveling waves and even the existence of a spreading speed for models with density-dependent dispersal probability. Some preliminary results and caveats were obtained in \cite{lut}. These questions remain the subject of our future investigation.

\end{document}